\documentclass[reqno,10pt]{amsart}
\usepackage{xspace}
\usepackage[bookmarksnumbered,colorlinks]{hyperref}
\usepackage{graphics}
\usepackage{fourier}
\newcommand{\bt}{\begin{theorem}}                     
\newcommand{\et}{\end{theorem}}                       
\newcommand{\bd}{\begin{definition}}                  
\newcommand{\ed}{\end{definition}}                    
\newcommand{\bl}{\begin{lemma}}                       
\newcommand{\el}{\end{lemma}}                                   
\newcommand{\bpr}{\begin{proposition}}                  
\newcommand{\epr}{\end{proposition}}                    
\newcommand{\bere}{\begin{remark}}                      
\newcommand{\ere}{\end{remark}}                         
\newcommand{\beq}{\begin{equation}}
\newcommand{\eeq}{\end{equation}}
\def\bal#1\eal{\begin{align}#1\end{align}}              
\def\baln#1\ealn{\begin{align*}#1\end{align*}}          
\def\bml#1\eml{\begin{multline}#1\end{multline}}        
\def\bmln#1\emln{\begin{multline*}#1\end{multline*}}  
\def\bga#1\ega{\begin{gather}#1\end{gather}}
\def\bgan#1\egan{\begin{gather*}#1\end{gather*}}
\newcommand{\de}{\mathrm{d}}                        
                     
\newcommand{\N}{\ensuremath{\mathbb{N}}\xspace}     
\newcommand{\R}{\ensuremath{\mathbb{R}}\xspace}     
\newcommand{\eps}{\varepsilon}                      
                   
\newcommand{\inte}{\int_0^1\!\!}

\newtheorem{theorem}{Theorem}[section]

\newtheorem{lemma}[theorem]{Lemma}
\newtheorem{proposition}[theorem]{Proposition}
\theoremstyle{definition}
\newtheorem{definition}[theorem]{Definition}
\theoremstyle{remark}
\newtheorem{remark}[theorem]{Remark}

\title[Addendum to Morse theory of causal geodesics]%
{Addendum to ``Morse theory of causal geodesics in a stationary spacetime via Morse theory of geodesics of a Finsler metric'', {A}nn. {I}nst. {H}. {P}oincar\'e {A}nal. {N}on  {L}in\'eaire, 27 (2010) 857--876}

\author[E. Caponio]{Erasmo Caponio}
\address{Dipartimento di Meccanica, Matematica e Managment. \hfill \break\indent
Politecnico di Bari, Via Orabona 4, 70125, Bari, ITALY}
\email{caponio@poliba.it}

\author[M. A. Javaloyes]{Miguel \'Angel Javaloyes}
\address{Departamento de Matem\'aticas. \hfill \break\indent 
Universidad de Murcia, Campus de Espinardo,
30100 Espinardo, Murcia, Spain}
\email{majava@um.es}

\author[A. Masiello]{Antonio Masiello}
\address{Dipartimento di Meccanica, Matematica e Managment. \hfill \break\indent
Politecnico di Bari, Via Orabona 4, 70125, Bari, ITALY}
\email{masiello@poliba.it}

\date{}

\begin{document}

\begin{abstract}
We give the details of the proof of equality (29) in \cite{CaJaMa09}
\end{abstract}
\maketitle
\section{Introduction}
In  \cite[Eq. (29)]{CaJaMa09}, we claim  that the relative homology groups 
$H_*(\tilde E^c_{|X}\cap \tilde O^*\!,\tilde E^c_{|X}\cap \tilde O^*\setminus\{0\})$ and $H_*(\tilde E^c\cap \tilde O^*\!,\tilde E^c\cap \tilde O^*\setminus\{0\})$ are isomorphic,   where, we recall, $X=C^1_0([0,1], U)$, $U$ is a neighbourhood of $0\in\R^n$, $\tilde E\colon H^1_0([0,1],U)\to\R$, $\tilde E(x)=\inte\tilde G(s,x,\dot x)\de s$,   $0\in  H^1_0([0,1],U)$ is a non-degenerate critical point of $\tilde E$,   $c=\tilde E(0)$, $\tilde E^c=\{x\in  H^1_0([0,1],U)\ |\ \tilde E(x)\leq c\}$ and $\tilde O^*$ is a neighbourhood of $0$ in   $H^1_0([0,1],U)$.  

For this we refer to the following  result by Palais \cite[Theorem 16]{Palais66a}:
\begin{theorem}[Palais, \cite{Palais66a}] 
Let $V_1$ and $V_2$ be two locally convex topological vector spaces, $f$ be a continuous linear map 
from $V_1$ onto a  dense linear subspace of $V_2$ and let  $O$ be an open subset of $V_2$ and $\tilde O=f^{-1}(O)$. If $V_1$ and $V_2$ are metrizable  then $\tilde f= f_{|\tilde O}\colon \tilde O\to O$ is a homotopy equivalence.
\end{theorem}
As a consequence, if $E$ is a Banach space which is dense and continuously immersed  in a Hilbert space $H$ and $(A,B)$ is a pair of open subsets of $H$ with $B\subset A$, then the relative homology groups $H_*(A,B)$ and $H_*(\tilde A,\tilde B)$, where $\tilde A= A\cap E$ and $\tilde B= B\cap E$, are isomorphic. 

In this  addendum  we would like to make clear how the above   result can   be applied to get  
\[
H_*(\tilde E^c_{|X}\cap \tilde O^*\!,\tilde E^c_{|X}\cap \tilde O^*\setminus\{0\}) \cong H_*(\tilde E^c\cap \tilde O^*\!,\tilde E^c\cap \tilde O^*\setminus\{0\}).
\]
Although it is not difficult to find some open subsets which are homotopically equivalent, with respect to the $H^1$ topology,  to the ones involved in the computations of the critical groups (cf., for example, \cite[Ch. III, Corollary 1.2]{Chang93}), it is not trivial to ensure, after applying Palais' result, that the intersections of these  subsets with $X$ continue to be homotopically equivalent
in the $C^1$ topology.  

Actually, the equality   between the critical groups   of a Dirichlet functional    with respect to the $H^1$ and $C^1$ topology is not a novelty   (cf. \cite{Chang93,Chang94,LiLiLi05}).  Anyway, there are some  issues for the    functional $\tilde E$    that we would like to   point out.     First,   $\tilde E$ is   not $C^2$ with respect to the $H^1$ topology (this is a very general phenomenon for smooth, at most quadratic  in the velocities Lagrangians cf. \cite[Prop. 3.2]{AbbSch09});   moreover, as   $\tilde G$ is  not everywhere twice differentiable,      $\tilde E$ is    also   not twice Gateaux differentiable at any {\em   non-$\tilde G$-regular  } curve (see Definition~\ref{nonreg}). Secondly,     although its  flow is well defined on $X$, the gradient of $\tilde E$ is not of the type identity plus a compact operator,  thus we cannot immediately state that it  possesses the   {\em retractible property} in  \cite[\S III]{Chang83}, which ensures that the deformation retracts  involved in 
the 
computation  of the critical groups are also continuous in $X$,  where the 
Palais-Smale condition does not hold.      To overcome this problem, we extend   a result in \cite{AbbSch09},   constructing   a    smooth   vector field, which is a pseudo-gradient in $\mathcal U\setminus \bar B(0,r)$, where $\mathcal U$ is a neighbourhood of $0$ in  
$H^1_
0([0,1],U)$ and $\bar B(0,r)$ is the closure of a ball, and whose   flow   satisfies  the retractible property.      

The proof we give in the next section (without Lemma~\ref{finitenumber}, which becomes superfluous) also holds for any {\em smooth} Lagrangian on  $[0,1]\times TM$, where $M$ is a finite dimensional manifold,  which is  fiberwise strongly  convex and has at most quadratic growth in each fibre. We can also consider, with minor modifications, more general boundary conditions as the curves joining two given submanifolds in $M$.  The Lagrangian action functional will be then defined  on  the Hilbert manifold of the $H^1$ curves between the two submanifolds.   As we have already mentioned above, such functional is in general not $C^2$. Assuming that  at least one of the submanifolds is compact and that all the critical points are non-degenerate,  we can obtain, as in \cite[Theorem 9]{CaJaMa09}, the Morse relations for the solutions of the corresponding Lagrangian system. In this case,  the number of the conjugate instants along a geodesic, counted with their multiplicity, is replaced by the number of the 
``focal instants'' with respect to one of the two submanifold (counted  with multiplicities) along a solution plus the index of a 
bilinear symmetric form related to the other submanifold, \cite{Crampin}. We recall that a Morse complex for the action functional of such 
kind of Lagrangian, whose homology is isomorphic to  the singular homology of the path space between the two submanifolds, has been obtained in 
\cite{AbbSch09}.
\section{Proof of the isomorphism between the critical groups in $H^1$ and $C^1$}
We recall that the Lagrangian 
$\tilde G\colon [0,1]\times U\times   \R  ^n\to [0,+\infty)$ is given by    
\[\tilde G(t,q,y)=F^2(\varphi(t,q),\de \varphi(t,q)[(1,y)]),\]  
  where $F$ is a Finsler metric on the $n$-dimensional smooth manifold $M$ and    $\varphi\colon [0,1] \times U\to M$ is defined as  
$\varphi(t,q)=\exp_{  \gamma_0(t) }P_t(q)$;   here,   $\exp$ is the exponential map with respect to any auxiliary Riemannian metric $h$ on $M$,  $\gamma_0$ is the geodesic of    $(M,F)$     in which    we want to compute the critical groups,    $P_t\colon U\to T_{  \gamma_0(t) } M$ is given by $P_t(q_1,\ldots,q_n)=\sum_{i=1}^{n}q_iE_i(t)$, where $\{E_i\}_{i\in\{1,\ldots,n\}}$     are $n$-orthonormal smooth vector fields along $\gamma_0 $ and $U$ is the   Euclidean    ball of radius $\rho/2$, where $\rho$ is the minimum of the injectivity radii (with respect to the metric $h$) at the points $  \gamma(t) $, $t\in[0,1]$. 
   
The set  $Z$ where $\tilde G$ is not twice differentiable is defined    by  the equation $ \de \varphi (t,q)[1,y]=0 $ and then it corresponds to the subset of $[0,1]\times U\times   \R^n$ where the Lagrangian $\tilde G(t, q, y)=F^2(\varphi (t,q),\de \varphi(t,q)[(1,y)]$ vanishes.
We recall also that for each $(t,q)\in[0,1]\times U$ there is only  one $y\in \R^n$ such that
$\de \varphi(t,q)[(1,y)]=0$. Indeed, $\de \varphi(t,q)[(1,y)]=\partial_t \varphi(t,q)+\partial_q\varphi(t,q)[y]$ and,  as $\partial_q\varphi(t,q)$ is   one-to-one , $y\in \R^n$ is the only vector such  that
\[\partial_q\varphi(t,q)[y]=-\partial_t \varphi(t,q).\]  

We recall also that the map $\varphi$ defines a smooth injective map $\varphi_*\colon H^1_0([0,1], U)\to \Omega_{p_0,q_0}(M)$, $\varphi_*(x)(t)=\varphi(t,x(t))$, such that $\tilde E=E\circ \varphi_*$, where $E$ is the energy functional of $F$, i.e. $E(\gamma)=\frac12 \inte F^2(\gamma,\dot\gamma)\de t$  and $\Omega_{p_0.q_0}$ is the Hilbert manifold of the $H^1$ curves on $M$ between $p_0$ and $q_0$.   Observe that the curve of constant value $0$ is mapped by $\varphi_*$ to the geodesic $\gamma_0$ (hence $0$ is a critical point of $\tilde E$).  

From the fact that $F^2$ is fiberwise positively homogeneous of degree $2$ and $\varphi$ is a smooth map, it follows that there exists a constant $c_1$, depending only on $U$, such that 
\beq
\|\tilde G_{qq}(s,q,y)\|\leq c_1(1+|y|^2),\quad\quad  \|\tilde G_{qy}(s,q,y)\|\leq c_1(1+|y|),\quad\quad  \|\tilde G_{yy}(s,q,y)\|\leq c_1,\label{esti}
\eeq 
for every $(s,q,y)\in [0,1]\times U\times \R^n\setminus Z$,    where $|\cdot|$ and $\|\cdot\|$ are, respectively,  the euclidean norm and the norm of bilinear forms on $\R^n$.

Moreover, since $F^2$ is fiberwise strongly convex, there exists a positive constant $c_2$ such that  
\beq \tilde G_{yy}(s,q,y) [w,w] \geq c_2  |w|^2, \label{esti2}\eeq
for each $(s,q,y)\in [0,1]\times U\times \R^n\setminus Z$   and $w\in \R^n$.  
\begin{definition}\label{nonreg}
A  curve $x\in H^1_0([0,1], U)$ is said {\em   $\tilde G$-regular } if  the set of points $  t \in[0,1]$ where $(  t ,x(  t ),\dot x(  t ))\in Z$ is negligible.
\end{definition}
Let $\alpha\colon \R^n\to \R$ be a smooth function such that $\alpha_{|U'}=1$, $\alpha_{|{U}^C}=0$, where $U'$ is an open subset of $\R^n$ such that $0\in U'$ and $\bar U'\subset U$. Consider the Lagrangian $\mathcal L\colon\R\times\R^n\times \R^n\to \R$, $\mathcal L(t,q,y)=\alpha(q)\tilde G(t,q,y)+(1-\alpha(q))|y|^2$.
Clearly,   by the definition of $\alpha$,   $0$ is   also a critical point   of   the action functional $\mathcal A_{\mathcal L}(x)=\frac{1}{2}\inte \mathcal L(s,x,\dot x)\de s$. Notice also that, like $\tilde E$, 
 $\mathcal A_{\mathcal L}\colon H^1_0([0,1],\R^n)\to\R$ is a $C^1$ functional with locally Lipschitz differential.

Let  $ \mathcal B$   be    a closed ball in $H^1_0([0,1],\R^n)$,   centred in $0$ and   containing  curves that have support in $U'$.

As   $\mathcal L=\tilde G$ on $\R\times U'\times \R^n$, we have that $\mathcal A_{\mathcal L}|_{\mathcal B}=\tilde E|_{\mathcal B}$.  Since $\tilde E$ satisfies the Palais-Smale condition (see \cite{CaJaMa07a}), we also have that $\mathcal A_{\mathcal L}$ satisfies the Palais-Smale condition in $ \mathcal B $.

  Moreover, from \eqref{esti} it follows 
that  $\tilde E$ is twice Gateaux differentiable at any    $\tilde G$-regular   curve $x\in H^1_0([0,1],U)$  and then  the same property is  satisfied by
$\mathcal A_{\mathcal L}$.

Observe that, as the endpoints of the geodesic $\gamma_0$ are not conjugate, then we can assume that $\mathcal B$ is an isolating neighbourhood of the critical point $0$.  
Moreover, the non-conjugacy assumption  
implies also that $0$ is a non-degenerate critical point of $\tilde E$,   that is,  the kernel of the operator $A$,  which represents the second Gateaux differential at $0$ of both $\tilde E$ and $\mathcal A_{\mathcal L}$, with respect to the scalar product $\langle \cdot, \cdot\rangle$ in $H^1_0([0,1],   \R ^n)$, is empty.

The following proposition has been obtained in \cite[Lemma 4.1 and formula (4.8)]{AbbSch09} for the action functional of a   $C^2$, time-dependent,  fiberwise strongly convex, at most quadratic in the velocities,  Lagrangian on $TM$.     
\bpr\label{abbo}
There exist a neighbourhood $\mathcal U'$ of $\,\, 0$ in  $H^1_0([0,1],\R^n)$ (that we can assume it is contained in $\mathcal B$) and a positive constant $\mu_0$, such that the linear vector field $x\in\mathcal U'\mapsto Ax$,  satisfies  the inequality
\beq \de \mathcal A_{\mathcal L}(x)[Ax]\geq \mu_0\|\nabla\mathcal A_{\mathcal L}(x)\|^2_0,\label{2pg}\eeq 
for each $x\in \mathcal U'$.
\epr
Here $\|\cdot\|_0$ is the $H^1_0$ norm.  
In our setting, the Lagrangian $\mathcal L$ is  not twice differentiable on $Z\subset TM$ and this leads to some differences between the proof of \cite[Lemma 4.1]{AbbSch09} and ours , which   we  outline  in lemmata \ref{finitenumber},  \ref{claim2} and \ref{claim3}. 

\bl\label{finitenumber}
Let $x$ be a  smooth curve (non necessarily   $\tilde G$-regular  ) in $H^1_0([0,1],U)$. Then  the curves $t\in[0,1]\mapsto  s x(t)$ can be   non-$\tilde G$-regular   for $s$ in a subset of $[0,1]$ which is at most countable.  
\el
\begin{proof} 
We recall    that  $w \in H^1_0([0,1], U)$,   $w=w(t)$, is not   $\tilde G$-regular   if $(t, w(t), \dot w(t))\in Z$ for each $t$ in a subset of positive  Lebesgue measure in $[0,1]$.
Now, for $x\colon [0,1]\to U$,  smooth and $x(0)=x(1)=0$,
let us consider the map $f\colon [0,1]\times [0,1]\to M$ defined as
$f(s,t)=\varphi(t,sx(t))$. Observe that for each $\bar t\in [0,1]$, $s\mapsto f(s, \bar t)$ is the affinely parametrized geodesic $\sigma_{\bar t}$
of the Riemannian metric $h$ defined by $\sigma_{\bar t}(s)=\varphi(\bar t,sx(\bar t))=\exp_{\gamma(\bar t)}(s x(\bar t))$ (for $\bar t=0$ and $\bar t=1$ the geodesics are constant)  while, for each $\bar s\in [0,1]$, $t\mapsto f(\bar s, t)$ is the curve $\gamma_{\bar s}$ corresponding to $\bar s x$ by the map  $\varphi_*$ (for $\bar s=0$ and $\bar s=1$, we get respectively   $\gamma_0$,  the geodesic of $(M, F)$,   and the curve $\gamma_1=
\varphi_*(x)$). 
Thus $f=f(s,t)$ defines a geodesic congruence and, then, $s\mapsto J_t(s)=\partial_tf(s,t)=\dot \gamma_s(t)$ defines a Jacobi field along $\sigma_t$ for each $t\in (0,1)$ where $x(t)\neq 0$. Observe that at the instants $\bar t$ where $x(\bar t)=0$ (if they exist), $\sigma_{\bar t}$ is constant and equal to $\gamma_0(\bar t)$. Since there is only one $y\in   \R^n$ such that $(\bar t, 0, y)\in Z$ and such $y$ cannot be equal to $0$ (otherwise $0=\de \varphi (\bar t,0)[1,0]=\partial_t\varphi (\bar t,0)+\partial_q\varphi(\bar t,0)[0]=\partial_t\varphi (\bar t,0)=\dot\gamma_0(\bar t)\neq 0$),  there can be at most one $s\in (0,1]$ such that $(\bar t, sx(\bar t), s\dot x(\bar t))=(\bar t, 0, s\dot x(\bar t))\in Z$. Now let us assume that for $s, s'\in (0,1]$, $s\neq s'$, the curves $sx$ and $s'x$ are not   $\tilde G$-regular.   From what we have recalled above, this is equivalent to the fact that the curve $\gamma_s$ and $\gamma_{s'}$ have velocity vector fields vanishing on, respectively, $Z_s\subset [0,1]
$ and $Z_{s'}\subset [0,1]$ with $|Z_s|, |Z_{s'}|>0$. We claim that $Z_s\cap Z_{s'}=\emptyset$. Indeed, if there exists $\bar t\in Z_s\cap Z_{s'}$,
then $x(\bar t)$ must be 
different from $0$ and this implies  that the Jacobi field $J_{\bar t}$ is well defined and equal to $0$ at the instants $s$ and $s'$. Thus the points $\sigma_{\bar t}(s)$ and $\sigma_{\bar t}(s')$ are conjugate along $\sigma_{\bar t}$, but this is impossible (see, e.g.,  \cite[Prop. 2.2,  p. 267]{Carmo92}) because such geodesic has length less than the injectivity radius at $\gamma_0(\bar t)$. Therefore  
the set $\mathcal Z$  of $s\in[0,1]$ such that $|Z_s|>0$  is at most countable.  
Indeed, by contradiction, assume that $\mathcal Z$ is uncountable and consider the set
$A_h=\{s\in (0,1]  : |Z_s|>\frac{1}{h}\}$. Since $\cup_{h\in \N} A_h=\mathcal Z$,  there must exist at least one $k\in \N$ such that $A_k$  is     uncountable. Thus, for infinitely many $s\in [0,1]$, we would have disjoint subsets $Z_s\subset [0,1]$ having measure greater than $\frac{1}{h}$,   which is impossible. 
\end{proof}
\bere\label{regdense}  
From  Lemma~\ref{finitenumber} it also follows that any smooth   non-$\tilde G$-regular    curve $x\in H^1_0([0,1], U)$ is the limit,  in the $H^1$ topology, of some sequence $(x_k)\subset H^1_0([0,1], U)$ of smooth   $\tilde G$-regular     curves. Indeed, it is enough to consider a sequence $(s_n)\subset [0,1]$ such that $s_n\to 1$  and $s_n x$ is   $\tilde G$-regular.   
\ere
\bere\label{varie}  From  \eqref{esti2}, the second Gateaux differential of $\tilde E$ at a   $\tilde G$-regular   curve $x$ is represented by a linear bounded self-adjoint operator on $H^1_0([0,1],\R^n)$ of the type $A_x=B_x+ K_x$ where $B_x$ is a strictly positive definite operator and $K_x$ is compact. Moreover from \eqref{esti}, if  a sequence of   $\tilde G$-regular   curves  $\{x_n\}$ converges to a   $\tilde G$-regular   curve $x$ in the $H^1$ topology then $K_{x_n}$ converges to  $K_x$ in the norm topology of the bounded operators  and $B_{x_n}$ converges strongly to $B_x$, i.e. $B_{x_n}[\xi]\to B_x[\xi]$ for each $\xi\in H^1_0([0,1],\R^n)$ (cf. claim 1 and 2 of the proof of Lemma~4.1 in \cite{AbbSch09}). We recall that from \cite[Lemma 2]{CaJaMa09}, $A\equiv A_0$ is given by $I+K$ (that is, $B_0$ is the identity operator).  
\ere
  
The following two results are the analogous of,  respectively, Eq. (4.5) and Claim 3 in \cite{AbbSch09}. 
 
\bl\label{claim2}
  Let $(x_n)\subset H^1_0([0,1], U)$ be a sequence of smooth   $\tilde G$-regular     curves   such that $x_n\to 0$ in  the $H^1$ topology. Then 
\[
\de \tilde E(x_n)[Ax_n]=\inte \left \langle (B_{sx_n}^{1/2}+ K_{sx_n})^2 x_n,x_n\right \rangle   \de s +o(\|x_n\|_0^2), \quad\text{as $n\to\infty$.}
\]
 
\el
\begin{proof}
  Eqs. \eqref{esti}-\eqref{esti2} imply that  $\tilde G(t,q,y)$ satisfies assumptions $(L1')$   and   $(L2')$ at page 605 of \cite{AbbSch09}, for each $(t,q,y)\in [0,1]\times U\times \R^n\setminus Z$.      
Hence the lemma follows arguing as in  \cite[Lemma~4.1]{AbbSch09}, taking into account that
\bal\label{inta}
\de \tilde E(x)[Ax]&=\left\langle \nabla \tilde E(x),x+K(x)\right\rangle
=\left(\inte \frac{\de}{\de s} \left\langle \nabla \tilde E(sx),x+K(x)\right\rangle\de s\right)\nonumber\\
&=\inte \left\langle (B_{sx}+K_{sx})x,x+K(x)\right\rangle\de s.
\eal 
In fact,
$\frac{\de}{\de s} \nabla \tilde E(s x) = (B_{sx}+ K_{sx})[x]$   
at the points $s$ where the curve $t\in[0,1]\mapsto sx(t)$ is   $\tilde G$-regular.  
From Lemma~\ref{finitenumber},    the set of points  $s\in[0,1]$ where $sx$ is   not $\tilde G$-regular   is at most countable.  \end{proof}
  The next lemma follows as in Claim 3 of \cite[Lemma~4.1]{AbbSch09}, recalling Remark~\ref{varie} and the fact that $0$ is a non-degenerate critical point of $\tilde E$. 
\bl\label{claim3}
  There exist a number  $\mu>0$   and a neighbourhood $\mathcal U''$ of $0$ in $H^1_0([0,1], U)$ such that, for each smooth and   $\tilde G$-regular   curve $x\in \mathcal U''$, the spectrum of the self-adjoint    operator   $B_x^{1/2}+K_x$ is disjoint from $[-\mu, \mu]$.
  \el
\begin{proof}[Proof of Proposition~\ref{abbo}]
  Since $\mathcal A_{\mathcal L}|_{\mathcal B}=\tilde E|_{\mathcal B}$, it is enough to   prove   the proposition for the functional $\tilde E$. From Lemmata~\ref{claim2} and \ref{claim3}, we get that there exists 
a positive constant $\mu_1$, such that 
\beq \de \tilde E(x)[A x]\geq \mu_1\|x\|_0^2,\label{2pg1}\eeq 
for   each  smooth   $\tilde G$-regular   curve    $x\in \mathcal U''$. From Remark~\ref{regdense} and the continuity of $\de \tilde E$ and $A$ with respect to the $H^1$ topology, inequality \eqref{2pg1} can be extended
to any smooth curve in $\mathcal U''$ and then, since smooth curves are dense in $H^1_0([0,1],U)$, to any $x\in \mathcal U''$.
As $\nabla \tilde E$ is a locally Lipschitz field and $\nabla \tilde E(0)=0$, we get
\[\de \tilde E(x)[A x]\geq \mu_0 \|\nabla \tilde E(x)\|_0^2,\]
for some positive constant $\mu_0$ and for all $x$ in some neighbourhood $\mathcal U'$ of $0$. 
 
\end{proof}
Now let $\eta_0\colon H^1_0([0,1],\R^n)\to [0,1]$ be a smooth bump function such that $\mathrm{supp}\, \eta_0\subset \mathcal U'$ and $\eta_0(x)=1$, for all $x\in\mathcal U$, where $\mathcal U$ is an open neighbourhood of $0$ in $H^1_0([0,1],\R^n)$ with $\overline{\mathcal U}\subset \mathcal U'$. Let us consider the vector field on 
$H^1_0([0,1],\R^n)$ defined as 
\[Y(x)=-\eta_0(x) Ax-\big(1-\eta_0(x)\big)\nabla\mathcal A_{\mathcal L}(x).\]
We point out that we cannot state that $Y$ is a pseudo-gradient vector field
because we are not able to prove that 
\beq\label{1pg}
\|Ax\|_0\leq  \mu_2   \| \de A_{\mathcal L}(x) \|_0,
\eeq
for some constant $  \mu_2  >\mu_0$ and all $x$ in some neighbourhood of $0$.\footnote{\label{corona}Actually using that $\mathcal A_{\mathcal L}$ satisfies the Palais-Smale condition and $0$ is an isolated critical point of $\mathcal A_{\mathcal L}$, we can prove that $Y$ satisfies \eqref{1pg} in any open subset $\mathcal U\setminus \bar B(0,r)$, where $B(0,r)$ is an open ball strictly contained  in $\mathcal U$, for a constant   $\mu_2$   depending on  $\mathcal U\setminus \bar B(0,r)$.}
Anyway  \eqref{2pg} implies  that  $Y$  satisfies the inequality
\beq
\de \mathcal A_{\mathcal L}(x)[Y(x)]\leq -\mu\|\nabla\mathcal A_{\mathcal L}(x)\|_0^2,\label{2pg2}\eeq 
for each $x\in H^1_0([0,1],\R^n)$, where $\mu=\min\{\mu_0,1\}$.   As we will show in Lemma~\ref{mw8.1}, inequality   \eqref{2pg2} (together with the remark in footnote~\ref{corona}) is enough to get a deformation result as in \cite[Lemma 8.3]{MawWil89}. 
For all $x\in H^1_0([0,1],\R^n)$, let $(\omega^-(x), \omega^+(x))$ be the maximal interval of definition of the solution of 
\beq\label{flow}
\begin{cases}
\dot\psi=  Y(\psi) ,\\
\psi(0)=x.
  \end{cases}
\eeq
 Observe that this problem is well-defined because $Y$ is a locally Lipschitz vector field in $H^1_0([0,1],\R^n)$, since $A$ and $\nabla\mathcal A_{\mathcal L}$ are. Furthermore,  \eqref{2pg2} implies that $\mathcal A_{\mathcal L}$ is decreasing along the flow of  $Y$  and as, $ Y_{|\mathcal U} =-A  =-I-K  $,   such flow   is given by
\begin{equation}\label{flowK}
\psi(x,t)=e^{-t}x-\int_0^te^{-t+s}K(\psi(x,s))\de s
\end{equation}
for $x\in \mathcal U$, whereas $\psi(x,t)\in  \mathcal U$.
The following lemma is an adaptation   
of Lemma 8.1 in \cite{MawWil89} to the flow of the vector field $Y$. 
\bl\label{mw8.1}
Let $\mathcal V$ be a closed neighbourhood of $0$ contained in $\mathcal U$. Then there exist $\eps>0$ and an open neighbourhood $O'\subset \mathcal V$ of 
$0$ in  $H^1_0([0,1],\R^n)$ such that  if $x\in O'$, then the solution  $\psi (x, \cdot )$ of \eqref{flow} either stays in $\mathcal V$ for all $t\in [0, +\infty)$ or it stays in $\mathcal V$ at least until $\mathcal A_{\mathcal L}(\psi (x, t))$ becomes less than $c-\eps$, (where $c=\mathcal A_{\mathcal L}(0)=\tilde E(0)$).
\el
 \begin{proof}
Observe that,   since   $Y_{|\mathcal V}=-A$,  $\psi(x,\cdot)$ is defined for all times until it lies in $\mathcal V$. Let $B(0,\rho)$ be the ball of radius $\rho$  centred  at $0$ such that $\bar B(0,\rho)\subset \mathcal{V}$ and let 
\[\mathcal C=\{x\in H^1_0([0,1],\R^n):\frac{\rho}{2}\leq\|x\|_0\leq \rho\}.\]
Since $\mathcal C\subset \mathcal B$, it is free of critical points and then
\begin{equation}\label{deltainf}
\delta=\inf_{x\in \mathcal C}\|\nabla \mathcal A_{\mathcal L}(x)\|_0>0,
\end{equation}
because $\mathcal A_{\mathcal L}$ satisfies the Palais-Smale condition on $\mathcal C$. Moreover
\beq\label{fromps} 
\|Y(x)\|_0=\|A x\|_0\leq \rho \|A\|_0 \leq \frac{\rho  \|A\|_0 }{\delta}\|\nabla\mathcal A_{\mathcal L}(x)\|_0,
\eeq
for each $x\in\mathcal C$. Let $\nu:= \frac{\rho  \|A\|_0 }{\delta}$ and
$O'=B(0,\rho/2)\cap \mathring{\widearc{\mathcal A_{\mathcal L}^{c+\frac{\mu\delta\rho}{4\nu}}}}$.  If $x\in O'$ is such that  $\psi(x,\bar t)$ does not belong to $\mathcal V$ for some $ \bar t>0 $, then there exist $0<t_1< t_2<\omega^+(x)$ such that
$\psi(x,t)\in \mathcal C$, for all $t\in (t_1,t_2)$ and $\|\psi(x,t_1)\|_0=\rho/2$, $\|\psi(x,t_2)\|_0=\rho$.  It follows that 
\bal
\mathcal A_{\mathcal L}(\psi(x,t_2))&=  \mathcal A_{\mathcal L}(\psi(x,t_1))+\int_{t_1}^{t_2}\de \mathcal A_{\mathcal L}(\psi(x,t))[Y(\psi(x,t))]\de t\nonumber\\
&\leq \mathcal A_{\mathcal L}(x)-\mu\int_{t_1}^{t_2}\|\nabla\mathcal A_{\mathcal L}(\psi(x,t))\|_0^2\de t\label{grad}\\
&\leq c+\frac{\mu\delta\rho}{4\nu}-\mu\delta\int_{t_1}^{t_2}\|\nabla\mathcal A_{\mathcal L}(\psi(x,t))\|_0\de t\nonumber\\
&\leq c+\frac{\mu\delta\rho}{4\nu} - \frac{\mu\delta}{\nu}\int_{t_1}^{t_2}\|Y(\psi(x,t))\|_0\de t\nonumber\\
&\leq c+\frac{\mu\delta\rho}{4\nu} -\frac{\mu\delta}{\nu}(\|\psi(x,t_2)\|_0-\|\psi(x,t_1)\|_0)\nonumber\\
&=c+\frac{\mu\delta\rho}{4\nu} -\frac{\mu\delta\rho}{2\nu}=c-\frac{\mu\delta\rho}{4\nu}.\nonumber
\eal
In the first inequality above,   we have used the fact that $\mathcal A_{\mathcal L}$ is decreasing in the flow of \eqref{flow} and inequality \eqref{2pg2};   in the second one,   the fact that $x\in O'\subseteq\mathcal A_{\mathcal L}^{c+\frac{\mu\delta\rho}{4\nu}}$ and \eqref{deltainf};     in the third one, inequality   \eqref{fromps};   in the last one, the following chain   of inequalities:
\[
\int_{t_1}^{t_2}\|Y(\psi(x,t))\|_0\de t=\int_{t_1}^{t_2}\|\dot\psi(x,t))\|_0\de t\geq \|\int_{t_1}^{t_2}\dot\psi(x,t)\de t\|_0\geq \|\psi(x,t_2)\|_0-\|\psi(x,t_1)\|_0.
\]
Thus the  conclusion  follows with $\eps= \frac{\mu\delta\rho}{4\nu} $.
\end{proof}
Let $V$ be  the subset of  $H^1_0([0,1],\R^n)$ given as $V=\displaystyle \bigcup_{x\in O'}\psi(x, [0,\omega^+(x))$,   where $O'$ is the neighbourhood of $0$ associated to $\mathcal V$ by Lemma~\ref{mw8.1}.   Since $O'$ is open, from standard  results  in ODE theory (cf. for example \cite[Corollary 4.2.10]{Lang95}), $V$ is also an open subset of $H^1_0([0,1],\R^n)$.  From  Lemma~\ref{mw8.1},  $\mathcal A_{\mathcal L}^{-1}\big((c-\eps,c+\eps)\big)\cap V\setminus \{0\}$ is contained in $\mathcal V\subset \mathcal U$ and it is free of critical points. 
\bl\label{mw8.2}
For every $x\in \mathcal A_{\mathcal L}^{-1}\big([c,c+\eps)\big)\cap V$, either there exists a unique $  T(x)  \in [0,\omega^+(x))$ such that $\mathcal A_{\mathcal L}\big (\psi(x,  T(x) )\big )=c$ or $\omega^+(x)=+\infty$ and 
$\psi(x,t)\to 0$, in $H^1_0([0,1],\R^n)$, as $t\to +\infty$.
\el
\begin{proof}
If $\mathcal A_{\mathcal L}(\psi(x,t))>c$, for all $t\in [0,\omega^+(x))$, then from Lemma~\ref{mw8.1},   $\omega^+(x)=+\infty$ and $\psi(x,t)\in \mathcal V$, for each  $t\in [0,+\infty)$.  
From inequality \eqref{grad},
\[\int_{0}^{+\infty}\|\nabla\mathcal A_{\mathcal L}(\psi(x,t))\|_0^2\de t\leq \frac{1}{\mu}\left(\mathcal A_{\mathcal L}(x)-c\right)<+\infty,\]
hence $\lim\inf_{t\to+\infty}\|\nabla\mathcal A_{\mathcal L}(\psi(x,t))\|_0^2=0$ and the Palais-Smale condition implies the existence of a sequence $\{t_n\}$ converging to $+\infty$ such that $\psi(x,t_n)\to 0$. Hence the  conclusion  follows from Lemma~\ref{mw8.1}.
\end{proof}
By Lemmata~\ref{mw8.1}  and \ref{mw8.2}, as in \cite[Lemma 8.3]{MawWil89}, we get  that  $\mathcal A_{\mathcal L}^{c}\cap V$ is a strong deformation retract of $\mathring{\widearc{\mathcal A_{\mathcal L}^{c+\eps/2}}}\cap V$. Analogously, $\mathcal A_{\mathcal L}^{c-\eps}\cap V$ is a strong deformation retract of both $\mathcal A_{\mathcal L}^{c}\cap V\setminus \{0\}$ and $\mathring{\widearc{\mathcal A_{\mathcal L}^{c-\eps/2}}}\cap V$.  Using that, for $A\subset B\subset C$,  if $B$ is  a strong deformation retract of $C$, then $H_*(B,A)\cong H_*(C,A)$ and if $A$ is a strong deformation retract of $B$, then  $H_*(C,A)\cong H_*(C,B)$ (for the last property,  see for example \cite[Property $H_6$--$\beta$]{Rothe73}), 
we obtain
\begin{equation}\label{intercrit}
H_*(\mathcal A_{\mathcal L}^{c}\cap V,\mathcal A_{\mathcal L}^{c}\cap V\setminus \{0\})\cong 
H_*(\mathring{\widearc{\mathcal A_{\mathcal L}^{c+\eps/2}}}\cap V,\mathring{\widearc{\mathcal A_{\mathcal L}^{c-\eps/2}}}\cap V).
\end{equation}
 
Let $O=\varphi_* (O')$ and $  \gamma_0=\varphi_*(0) $, then
\begin{align}
C_*(E,  \gamma_0 )&= H_*(E^c\cap O,E^c\cap O\setminus\{  \gamma_0  \})\cong H_*((E\circ\varphi_*)^c\cap O',(E\circ\varphi_*)^c\cap O'\setminus\{0\})\nonumber\\
&=H_*(\tilde E^c\cap O',\tilde E^c\cap O'\setminus\{0\})=H_*(\mathcal A_{\mathcal L}^c\cap O',\mathcal A_{\mathcal L}^c\cap O'\setminus\{0\})\nonumber\\
&\cong H_*(\mathcal A_{\mathcal L}^c\cap V,\mathcal A_{\mathcal L}^c\cap V\setminus\{0\}),\label{finalstep}
\end{align}
last  equivalence, by the excision property of the singular relative homology groups.
By Palais' theorem above we get 
\[H_*(\mathring{\widearc{\mathcal A_{\mathcal L}^{c+\eps/2}}}\cap V,\mathring{\widearc{\mathcal A_{\mathcal L}^{c-\eps/2}}}\cap V)\cong
H_*(\mathring{\widearc{\mathcal A_{\mathcal L}^{c+\eps/2}}}\mbox{|}_{C^1_0([0,1],\R^n)}\cap V,\mathring{\widearc{\mathcal A_{\mathcal L}^{c-\eps/2}}}\mbox{|}_{C^1_0([0,1],\R^n)}\cap V).\]
The above equivalence, together with \eqref{intercrit} and \eqref{finalstep}, implies that
\[C_*(E,  \gamma_0 )\cong H_*(\mathring{\widearc{\mathcal A_{\mathcal L}^{c+\eps/2}}}\mbox{|}_{C^1_0([0,1],\R^n)}\cap V,\mathring{\widearc{\mathcal A_{\mathcal L}^{c-\eps/2}}}\mbox{|}_{C^1_0([0,1],\R^n)}\cap V).\]

It remains to prove that these last relative homology groups are isomorphic to the critical groups in $X=C^1_0([0,1],U)$.
To this end, let us consider 
  the Cauchy problem \eqref{flow}, with $x\in C^1([0,1],\R^n)\cap \mathcal A_{\mathcal L}^{-1}\big((c-\eps/2,c+\eps/2)\big)\cap V$.
Since   $\mathcal A_{\mathcal L}^{-1}\big((c-\eps/2,c+\eps/2)\big)\cap V\subset  \mathcal V\subset \mathcal U$, it holds \eqref{flowK} and
 the orbit $\psi(x,\cdot)$, defined by $x$, is also in $C^1_0([0,1],\R^n)$.

As a consequence,  the strong deformation retracts that we   have  considered  above    are well defined in $C^1_0([0,1],\R^n)\times[0,1]$ and by the continuity of the flow \eqref{flowK} with respect to the $C^1$ topology, we immediately deduce that they are also continuous at each point different from $(0,1)$. Clearly, the continuity at the point $(0,1)$ with respect to the product topology of   $C^1_0([0,1],\R^n)$, with the $C^1$ topology, and $\R$, with the standard one,    comes into play only for   the   deformation map 
$\eta\colon\mathring{\widearc{\mathcal A_{\mathcal L}^{c+\eps/2}}}\cap V\times [0,1]\to \mathring{\widearc{\mathcal A_{\mathcal L}^{c+\eps/2}}}\cap V$
of $\mathring{\widearc{\mathcal A_{\mathcal L}^{c+\eps/2}}}\cap V$ in $\mathcal A_{\mathcal L}^{c}\cap V$, 
  which is given by
\[\eta(x,t)= \begin{cases}
\rho\left(x, \frac{t}{t-1}\right)          
            &\text{if $t\in [0,1)$},\\
\displaystyle \lim_{s\to+\infty}\rho(x,s)&\text{if $t=1$},
\end{cases}
\]
where $\rho\colon \mathring{\widearc{\mathcal A_{\mathcal L}^{c+\eps/2}}}\cap V \times [0,+\infty)\to \mathring{\widearc{\mathcal A_{\mathcal L}^{c+\eps/2}}}\cap V$ is the map defined as follows:
if $\mathcal A_{\mathcal L}(x)> c$ and there exists $T(x)>0$  such that $\mathcal A_{\mathcal L}\big(\psi(x,T(x))\big)=c$, then 
\[\rho(x,t)=\begin{cases}
\psi(x,t)&\text{if $t\in [0, T(x)]$},\\
\psi(x,T(x))&\text{if $t\in (T(x),+\infty)$},             
            \end{cases}\]
if $\psi(x,t)\to c$ as $t\to +\infty$, then $\rho(x,t)=\psi(x,t)$ and if $\mathcal A_{\mathcal L}(x)\leq c$, then $\rho(x,t)=x$, for all $t\in[0,+\infty)$.
  Since the flow $\psi_1$ of the linear vector field $x\mapsto -A x= -Ix-Kx$ is given  by \eqref{flowK}
and $K$ is bounded from $H^1_0([0,1],\R^n)$ to $C^1_0([0,1],\R^n)$, we  have 
\[\|\int_0^te^{-t+s}K(\psi_1(x,s))\de s\|_{C^1}\leq e^{-t}\int_0^t e^{s}\|K(\psi_1(x,s))\|_{C^1}\de s\leq C e^{-t}\int_0^t e^{s} \|\psi_1(x,s)\|_{0} \de s.\]  
Thus, if   $\psi(x,t)\to 0$   in $H^1$, as $t\to+\infty$,    then, from {Lemmata \ref{mw8.1} and \ref{mw8.2}}, $\psi(x,t)=\psi_1(x,t)$. Hence,   for every $  \eps >0$, there exists $\bar t>0$ such that for all $t>\bar t$, $ \|\psi(x,t)\|_{0}< \eps  $ and then the last function in the above inequalities can be estimated, for $t>\bar t$, as
\bmln e^{-t}\int_0^t e^{s} \|\psi(x,s)\|_{0} \de s=e^{-t}\int_0^{\bar t} e^{s} \|\psi(x,s)\|_{0} \de s+e^{-t}\int_{\bar t}^t e^{s} \|\psi(x,s)\|_{0} \de s\\
\leq e^{-t}\int_0^{\bar t} e^{s} \|\psi(x,s)\|_{0} \de s+  \eps  (1-e^{-t}e^{\bar t}).
\emln
Thus $ \psi(x,t)\to 0 $ also with respect to the $C^1$ topology, giving the continuity of  the map $\eta$   at the point    $(0,1)$ also with respect to the product of    such a topology and the Euclidean one on the interval $[0,1]$. 
 
In conclusion we have that the following groups are isomorphic
\bmln
H_*(\mathring{\widearc{\mathcal A_{\mathcal L}^{ c+\eps/2}}}\mbox{|}_{C^1_0([0,1],\R^n)}\cap V,\mathring{\widearc{\mathcal A_{\mathcal L}^{ c-\eps/2}}}\mbox{|}_{C^1_0([0,1],\R^n)}\cap V)\cong\\
 \cong H_*({\mathcal A_{\mathcal L}^{c}}\mbox{|}_{C^1_0([0,1],\R^n)}\cap V,{\mathcal A_{\mathcal L}^{c}}\mbox{|}_{C^1_0([0,1],\R^n)}\cap V\setminus\{0\}).\emln
By excision, these last relative homology groups are isomorphic to 
$H_*({\mathcal A_{\mathcal L}^{c}}\mbox{|}_{C^1_0([0,1],\R^n)}\cap O',{\mathcal A_{\mathcal L}^{c}}\mbox{|}_{C^1_0([0,1],\R^n)}\cap O'\setminus\{0\})$ and then, since the curves in $O'$ have their support in $U$, to
$H_*(\tilde E^{c}\mbox{|}_{C^1_0([0,1],U)}\cap O',\tilde E^{c}\mbox{|}_{C^1_0([0,1],U)}\cap O'\setminus\{0\})$. 


\begin{thebibliography}{10}
\bibitem{AbbSch09}
{\sc A.~Abbondandolo and M.~Schwarz}, {\em A smooth pseudo-gradient for the Lagrangian  action functional},
Advanced Nonlinear Studies, 9 (2009), pp.~597--623.

\bibitem{CaJaMa07a}
{\sc E.~Caponio, M.~A. Javaloyes, and A.~Masiello}, {\em On the energy
  functional on {F}insler manifolds and applications to stationary spacetimes}, 
Math. Ann., 351 (2011), pp.~365--392.

\bibitem{CaJaMa09}
{\sc E.~Caponio, M.~A. Javaloyes, and A.~Masiello}, {\em Morse theory of causal
  geodesics in a stationary spacetime via {M}orse theory of geodesics of a
  {F}insler metric}, {A}nn. {I}nst. {H}. {P}oincar\'e {A}nal. {N}on
  {L}in\'eaire, 27 (2010), pp.~857--876. 

\bibitem{Chang83}
  {\sc K.-C. Chang},
  {\em A variant mountain pass lemma},
  Sci. Sinica Ser. A, 26 (1983), pp.~1241--1255.



 
\bibitem{Chang93}
   {\sc K.-C. Chang},
  {\em Infinite-Dimensional {M}orse Theory and Multiple Solution Problems},
  {Birkh{\"a}user},Boston, MA, 1993.  


\bibitem{Chang94}
{\sc K.-C. Chang},
{\em $H\sp 1$ versus $C\sp 1$ isolated critical points},
{C. R. Acad. Sci. Paris S{\'e}r. I Math.},
319 (1994), pp.~441--446.
  
\bibitem{Crampin}
{\sc M. Crampin},
{\em The {M}orse index theorem for general end conditions},
Houston J. Math.,  27 (2001), pp.~807--821.
  
\bibitem{Carmo92}
{\sc M. P. do Carmo},
{\em Riemannian Geometry},
Birkh{\"a}user,Boston, MA, 1992.
  



\bibitem{Lang95}
{\sc S. Lang},
{\em {D}ifferential and {R}iemannian manifolds}, Graduate Texts in Mathematics, Springer-Verlag, New York, 1995.


\bibitem{LiLiLi05}
  {\sc C. Li, S.  Li and J. Liu},
  {\em Splitting theorem, {P}oincar\'e-{H}opf theorem and jumping nonlinear problems},
  J. Funct. Anal., 221 (2005), pp.~439--455.


\bibitem{MawWil89}
{\sc J.~Mawhin and M.~Willem}, {\em Critical point theory and {H}amiltonian
  systems},  Applied Mathematical Sciences, Springer-Verlag, New
  York, 1989.


\bibitem{Palais66a}
  {\sc R. S. Palais},
  {\em Homotopy theory of infinite dimensional manifolds},
  {Topology}, 5 (1966),
  pp.~1--16.

\bibitem{Rothe73}
  {\sc E. H. Rothe},
  {\em {M}orse theory in {H}ilbert space},
  Rocky Mountain J. Math.,
  3 (1973), pp. 251--274.
\end{thebibliography}
\end{document}